\documentclass[11pt]{article}

\usepackage{amssymb}
\usepackage{amsmath}
\usepackage{amsthm} 
\usepackage{mathtools}
\usepackage{bm} 
\usepackage{extarrows}
\usepackage{indentfirst}
\usepackage{stmaryrd}
\usepackage{tikz}
\usetikzlibrary{arrows, graphs}
\usepackage{setspace}
\def\ex{\mbox{ex}}

\newtheorem{thm}{Theorem}[section]

\newtheorem{lem}[thm]{Lemma}
\newtheorem{prop}{Proposition}
\newtheorem{conj}[thm]{Conjecture}

\newtheorem{claim}{Claim}

\begin{document}

\title{Graphs without rainbow cliques of orders four and five 
\thanks{The work was supported by the National Natural Science Foundation of China (No. 12071453),
	and the National Key R and D Program of China (2020YFA0713100), and the Innovation Program for Quantum
	Science and Technology, China (2021ZD0302902).}
}
\author{Yue Ma$^a$, \quad Xinmin Hou$^{a,b}$\\
\small $^a$School of Mathematical Sciences\\
\small University of Science and Technology of China\\
\small Hefei, Anhui 230026, China.\\
\small $^{a,b}$ Key Laboratory of Wu Wen-Tsun Mathematics\\
\small University of Science and Technology of China\\
\small Hefei, Anhui 230026, China.
}

\date{}

\maketitle

\begin{abstract}
Let $\mathcal{G}_n^k=\{G_1,G_2,\ldots,G_k\}$ be a multiset of graphs on vertex set $[n]$ and let $F$ be a fixed graph with edge set $F=\{e_1, e_2,\ldots, e_m\}$ and $k\ge m$. We say ${\mathcal{G}_n^k}$ is rainbow $F$-free if there is no $\{i_1, i_2,\ldots, i_{m}\}\subseteq[k]$ satisfying $e_j\in G_{i_j}$ for every $j\in[m]$. 
Let $\ex_k(n,F)$ be the maximum $\sum_{i=1}^{k}|G_i|$ among all the rainbow $F$-free multisets ${\mathcal{G}_n^k}$.
Keevash, Saks, Sudakov, and Verstra\"ete (2004) determined the exact value of $\ex_k(n, K_r)$  when $n$ is sufficiently large and proposed the conjecture that the results remain true when $n\ge Cr^2$ for some constant $C$. Recently, Frankl (2022) confirmed the conjecture for $r=3$ and all possible values of $n$. In this paper, we determine the exact value of $\ex_k(n, K_r)$ for $n\ge r-1$ when $r=4$ and $5$, i.e.  the conjecture of Keevash, Saks, Sudakov, and Verstra\"ete is true for $r\in\{4,5\}$.
\end{abstract}

\section{Introduction}
Let $G=(V,E)$ be a graph with vertex set $V$ and edge set $E\subset\binom{V}{2}$. We write $G$ instead of $E(G)$ in this article. For any $v\in V(G)$, let $d_{G}(v)=|\{e\in G: v\in e\}|$ be the degree of $v$ in $G$ and $\delta(G)=\min_{v\in V(G)}d_{G}(v)$. We may simply write $d(v)$ when it causes no confusion. 
Let $K_r$ denote a complete graph on $r$ vertices.
For disjoint sets $U_1,U_2,...,U_r$ ($r\ge 2$), let $K[U_1,U_2,...,U_r]$ denote a complete $r$-partite graph with partition sets $U_1,...,U_r$.
When $U_1,U_2,...,U_r\subseteq V(G)$, denote $G[U_1,...,U_r]=K[U_1,...,U_r]\cap G$.

The Tur\'an graph $T_{r-1}(n)$ is a complete $(r-1)$-partite graph on $n$ vertices whose partition sets have sizes as equal as possible. 
Let $t_{r-1}(n)=|T_{r-1}(n)|$.
The famous Tur\'an Theorem~\cite{Tr} states that $|G|\le t_{r-1}(n)$ for an $n$-vertex $K_r$-free $G$. 
In 1966, Simonovits~\cite{S68} extended the Tur\'an Theorem by showing that $|G|\le t_{r-1}(n)$ for an $n$-vertex $F$-free $G$ for sufficiently large $n$, where $F$ is a graph with chromatic number $r$ and there is some edge $e$ such that $F-e$ has chromatic number $r-1$ (such a graph $F$ is called {\em $r$-critical}). 
There are fruitful results in the field of classical Tur\'an problems (one can see good surveys, for example~\cite{Caen94,Furedi91,Kee2011}).

In 2004, Keevash, Saks, Sudakov, and Verstra\"ete~\cite{KSSV04} introduced a rainbow version of the Tur\'an-type problem. Let $F$ be a fixed graph with $E(F)=\{e_1, e_2,\ldots, e_{|F|}\}$. For integers $k\ge |F|$ and $n$, let ${\mathcal{G}_n^k}=\{G_1, G_2, ..., G_k\}$ be a multiset of graphs on vertex set $[n]$. We say ${\mathcal{G}_n^k}$ is {\em rainbow $F$-free}, or {\em RB$F$-free} if there is no $\{i_1, i_2, ..., i_{|F|}\}\subseteq[k]$ satisfying $e_j\in G_{i_j}$ for every $j\in[|F|]$.
The {\em rainbow Tur\'an-type problems} look for the maximum of $\sum_{i=1}^{k}|G_i|$ or $\prod_{i=1}^{k}|G_i|$ among all the RB$F$-free  multisets ${\mathcal{G}_n^k}$. Let $|{\mathcal{G}_n^k}|=\sum_{i=1}^{k}|G_i|$. The {\em rainbow Tur\'an number} of a fixed graph $F$ is defined as 
$$\ex_k(n, F)=\{|\mathcal{G}_n^k| : \mathcal{G}_n^k \text{ is RB$F$-free}\}.$$ 
A multiset of graphs $\mathcal{G}_n^k$ with $|{\mathcal{G}_n^k}|=\ex_k(n, F)$ is called an {\em extremal system} of rainbow $F$.

Denote $t\mathcal{K}_n$ (or $t\overline{\mathcal{K}}_n$) as the multiset consisting of $t$ copies of complete graphs $K_n$ (or empty graphs $\overline{K}_n$).
Let $t{\mathcal{T}_{r-1}(n)}$ be the multiset consisting of $t$ copies of the Tur\'an graph $T_{r-1}(n)$.

For $k\ge h$ and an $r$-critical graph $H$ with $h$ edges, it is not difficult to see that both ${\left(h-1\right)\mathcal{K}_{n}\cup\left(k-h+1\right)\overline{\mathcal{K}}_n}$ and ${ k\mathcal{T}_{r-1}(n)}$ contain no rainbow $H$. It is reasonable to consider them as the extremal structures for $H$. In fact, 
Keevash, Saks, Sudakov, Verstra\"ete~\cite{KSSV04} proposed the following conjecture.
\begin{conj} [Keevash, Saks, Sudakov, Verstra\"ete, 2004]\label{CONJ: 1} Let $r\ge 3$ and $H$ be an $r$-critical graph with $h$ edges. Suppose
$k\ge h$, $n$ is sufficiently large. Then 
$$\ex_k(n, H)=\max\left\{k\cdot t_{r-1}(n), \left(h-1\right)\binom{n}2\right\}.$$
Moreover, ${\left(h-1\right)\mathcal{K}_{n}\cup\left(k-h+1\right)\overline{\mathcal{K}}_n}$ and ${ k\mathcal{T}_{r-1}(n)}$ are the only extremal structures when $n$ is sufficiently large. 
\end{conj}

This conjecture has been confirmed when $r=3$ or $H=K_r$ for $r\ge 3$ in~\cite{KSSV04}. 
\begin{thm}[Keevash, Saks, Sudakov, Verstra\"ete~\cite{KSSV04}]\label{THM:KSSV}
Let $r\ge 2$, $k\ge {r\choose 2}$, $n> 10^4r^{34}$. Then 
$$\ex_k(n,K_r)=\begin{cases}
	k\cdot t_r(n), \text{ for } k\ge\frac 12(r^2-1),\\
	\left(\binom{r}{2}-1\right)\binom{n}{2}, \text{ for } {r\choose 2}\le k <\frac 12(r^2-1).
	\end{cases}$$
Moreover, ${\left(h-1\right)\mathcal{K}_{n}+\left(k-h+1\right)\overline{\mathcal{K}}_n}$ and ${ k\mathcal{T}_{r-1}(n)}$ are the only extremal structures. 
\end{thm}

However, the lower bound $n> 10^4r^{34}$ is unnatural. Keevash, Saks, Sudakov, Verstra\"ete~\cite{KSSV04} also proposed the following interesting conjecture: Theorem~\ref{THM:KSSV} probably remains true even for $ n\ge Cr^2$ for some constant $C$. We may furthermore hope   Theorem~\ref{THM:KSSV} also holds for all possible values of $n$. 
\begin{conj}[Keevash, Saks, Sudakov, Verstra\"ete~\cite{KSSV04}]\label{CONJ: C2}
Let $n\ge r-1\ge 2$, $k\ge\binom{r}{2}$. 
Then 
$${\ex_k(n, K_r)}\le\max\left\{\left(\binom{r}{2}-1\right)\binom{n}{2}, k\cdot t_{r-1}(n)\right\}\mbox{.}$$
\end{conj}

Recently, Frankl~\cite{K3} proved Conjecture~\ref{CONJ: C2} for $r=3$ and all possible values of $n$. 
In this paper, we continue to show that Conjecture~\ref{CONJ: C2} holds for $r=4,5$.
\begin{thm}\label{THM: c2r45}
Let $r\in\{4,5\}$, $n\ge r-1$, and $k\ge\binom{r}{2}$. 
Then 
$${\ex_k(n, K_r)}\le\max\left\{\left(\binom{r}{2}-1\right)\binom{n}{2}, k\cdot t_{r-1}(n)\right\}\mbox{.}$$
\end{thm}

The rest of this article is arranged as follows. In Section 2, we transfer the RB$F$-free system ${\mathcal{G}_n^k}=\{G_1, G_2,\ldots, G_{k}\}$ into a weighted graph  and give a weighted version (Conjecture~\ref{wightmain}) of Conjecture~\ref{CONJ: C2} and show that Conjecture~\ref{wightmain} implies Conjecture~\ref{CONJ: C2}. In Section 3, we prove Conjecture~\ref{wightmain} when $r=4$ and $5$, which implies Theorem~\ref{THM: c2r45}. We provide some discussion  in the last section.

\section{The $k$-weighted graph and the $k$-weighted version of Conjecture~\ref{CONJ: C2}}
The following lemma can be found in~\cite{KSSV04}.
\begin{lem}[\cite{KSSV04}]\label{subset}
Suppose the multiset of graphs ${\mathcal{G}_n^k}=\{G_1, G_2,\ldots, G_{k}\}$ on vertex set $[n]$ is RB$F$-free. Then there exists another RB$F$-free multiset of graphs $\mathcal{H}_n^k=\{H_1, H_2,\ldots, H_{k}\}$ on vertex set $[n]$ satisfying:

(a) $\cup_{i=1}^kH_i=\cup_{i=1}^k G_i$;

(b) $H_k\subseteq H_{k-1}\subseteq\ldots\subseteq H_1$.
\end{lem}

For some integer $\ell$ and a graph $H$, an {\it $\ell$-weighting} on $H$ is a mapping $g: E(H)\to \{0, 1, 2, \ldots, \ell\}$. We call the system $\mathcal{H}=(H,g)$ an {\it $\ell$-weighted graph}. 
For a subgraph $H'\subseteq H$, let $g(H')=\sum_{e\in H'}g(e)$ be the {\it weight} of $H'$ in $\mathcal{H}$. Call $(a_1,a_2,...,a_{|H'|})$ in a nondecreasing order  a {\it weight sequence} of $H'$ if $H'=\{e_1, e_2, \ldots, e_{|H'|}\}$ and $g(e_i)=a_i$ for any $1\le i\le |H'|$. A sequence  $(b_1,b_2,...,b_{|H'|})$ in nondecreasing  order is called a {\it weight sequence bound} of $H'$ if  $(a_1,a_2,...,a_{|H'|})\ge (b_1,b_2,...,b_{|H'|})$ with lexicographical order for every weight sequence $(a_1,a_2,...,a_{|H'|})$ of $H'$.

We propose the following conjecture of a $k$-weighted graph that implies Conjecture~\ref{CONJ: C2}.
\begin{conj}\label{wightmain}
Let $n\ge r\ge 3$ and $k\ge\binom{r}{2}$. Let $\mathcal{G}=(G_0,f)$ be a $k$-weighted graph, where $G_0\cong K_n$ on vertex set $[n]$. 
Suppose that $G_0$ contains no complete subgraph $K_r$ with weight sequence bound $(1,2,...,\binom{r}{2})$. Let $k_2=\left\lceil\frac{\left(\binom{r}{2}-1\right)\binom{n}{2}}{t_{r-1}(n)}\right\rceil$ and $k_1=k_2-1$. Then the following statements hold.

(i) If $k=k_1$, then $f(G_0)\le(\binom{r}{2}-1)\binom{n}{2}$.

(ii) If $k=k_2$, then $f(G_0)\le k\cdot t_{r-1}(n)$.
\end{conj}

\begin{thm}\label{wightimpliesc2}
Conjecture~\ref{wightmain} implies Conjecture~\ref{CONJ: C2}.
\end{thm}
\begin{proof}
Let ${\mathcal{G}_n^k}=\{G_1, G_2,..., G_{k}\}$ be a RB$K_r$-free system with maximum $|\mathcal{G}_n^k|$.
By Lemma~\ref{subset}, we may assume $G_k\subseteq G_{k-1}\subseteq\ldots\subseteq G_1$. Let $G_0=K_n$ and define a $k$-weighting $f$ on $G_0$ by setting $f(e)=\max_{e\in G_i}i$ for every edge $e\in G_0$. Thus $f(e)$ is the number of graphs containing $e$ in $\mathcal{G}_n^k$  because $G_k\subseteq\ldots\subseteq G_1\subseteq G_0$. Therefore, we have $f(G_0)=|{\mathcal{G}_n^k}|$. We may assume $k\ge k_1$, otherwise, we can add $k_1-k$ copies of the empty graphs $K_n^c$ into $\mathcal{G}_n^k$ with the property RB$K_r$-free  and the $k$-weighted graph $G_0$ unchanged.

Now we show that there is no $K_r$ in $G_0$ with weight sequence bound $(1,2,...,\binom{r}{2})$. Otherwise, suppose there is a copy $H$ of $K_r$ in $G_0$ with weight sequence $(a_1,a_2,...,a_{h})\ge (1,2,...,\binom{r}{2})$, where $h=|H|=\binom{r}{2}$ and assume $H=\{e_1,e_2,\ldots,e_{h}\}$ with $f(e_i)=a_i$ for any $i\in[h]$. Since $(1,2,...,\binom{r}{2})\le (a_1,a_2,...,a_{h})$, we have $f(e_i)=a_i\ge i$ for every $i\in[h]$. Thus, by the definition of $f$ and $G_k\subseteq\ldots\subseteq G_1\subseteq G_0$, we have $e_i\in G_i$ for every $i\in[\binom{r}{2}]$. This leads to a rainbow $H\cong K_r$ in ${\mathcal{G}_n^k}$, a contradiction.

Now we suppose Conjecture~\ref{wightmain} is true.  If $k=k_1$ or $k_2(=k_1+1)$, then $|{\mathcal{G}_n^k}|=f(G_0)\le\min\left\{\left(\binom{r}{2}-1\right)\binom{n}{2}, k\cdot t_{r-1}(n)\right\}$. Thus Conjecture~\ref{CONJ: C2} holds when $k=k_1$ or $k_2$. 
If $k>k_2$, let ${\mathcal{G'}_n^{k_2}}=\{G_1, ..., G_{k_2}\}$. Obviously, ${\mathcal{G'}_n^{k_2}}$ is still RB$K_r$-free. By Conjecture~\ref{wightmain}, $|{\mathcal{G'}_n^{k_2}}|\le k_2\cdot t_{r-1}(n)$. Since $G_k \subseteq\ldots\subseteq G_{k_2} \subseteq\ldots\subseteq G_1$, we have $|G_i| \le |G_j|$ for every $k_2+1\le i\le k$ and $1\le j\le k_2$. Hence $$\sum_{i=k_2+1}^k|G_i|\le\frac{k-k_2}{k_2}\sum_{j=1}^{k_2}|G_j|=\frac{k-k_2}{k_2}|{\mathcal{G'}_n^{k_2}}|\le(k-k_2)t_{r-1}(n),$$ which implies that $$|{\mathcal{G}_n^k}|=\sum_{i=k_2+1}^k|G_i|+\sum_{j=1}^{k_2}|G_j|\le (k-k_2)t_{r-1}(n)+k_2\cdot t_{r-1}(n)=k\cdot t_{r-1}(n).$$ We are done.
\end{proof}
Therefore, to prove Theorem~\ref{THM: c2r45}, it is sufficient to show Conjecture~\ref{wightmain} for $r=4,5$.

\section{Conjecture~\ref{wightmain} holds when $r=4$ or $5$}
Before the proof, we first give some computational properties about $t_{r-1}(n)$.
\begin{prop}\label{bound}
For integers $n\ge r-1\ge 3$, let $k_2=\left\lceil\frac{\left(\binom{r}{2}-1\right)\binom{n}{2}}{t_{r-1}(n)}\right\rceil$. Then the following hold.

(i) If $n\ge r+1$, then $k_2\ge\binom{r}{2}+1$.

(ii) For $1\le s\le\max\{r-1,n-1\}$, $t_{r-1}(n)-t_{r-1}(n-s)\ge\binom{s}{2}+\frac{r-2}{r-1}s(n-s)$.

(iii) For $1\le s\le n-1$, $\frac{t_{r-1}(n-s)}{t_{r-1}(n)}\ge\frac{\binom{n-s}{2}}{\binom{n}{2}}$.

(iv) If $n\ge r+1$, for $1\le s\le r-1$, $$\left(\binom{r}{2}-1\right)\left(\binom{n}{2}-\binom{n-s}{2}\right)\ge\binom{r}{2}\left(\binom{s}{2}+\frac{r-2}{r-1}s(n-s)\right)+n-s.$$
\end{prop}
\begin{proof}
We write $n=k_0(r-1)+m$, where $k_0=\lfloor\frac{n}{r-1}\rfloor$ and $0\le m\le r-2$.
\begin{eqnarray*}
t_{r-1}(n) &=& \binom{n}{2}-m\binom{k_0+1}{2}-(r-1-m)\binom{k_0}{2}\\
&=&\binom{n}{2}-\frac{k_0}{2}(m(k_0+1)+(r-1-m)(k_0-1))\\
&=&\binom{n}{2}-\frac{k_0}{2}(k_0(r-1)+2m-(r-1))\\
&=&\binom{n}{2}-\frac{n-m}{2(r-1)}(n+m-(r-1))\mbox{.}
\end{eqnarray*}
Let $g_r(m,n)=\binom{n}{2}-\frac{n-m}{2(r-1)}(n+m-(r-1))$. By simple calculation, we have $$t_{r-1}(n)=g_r(m,n)\le g_r(0,n)=\binom{n}{2}-\frac{n(n-(r-1))}{2(r-1)}.$$ Note that 
$$k_2=\left\lceil\frac{\left(\binom{r}{2}-1\right)\binom{n}{2}}{t_{r-1}(n)}\right\rceil=\left\lceil\left(\binom{r}{2}-1\right)\left(\frac{\binom{n}{2}}{g_r(m,n)}-1\right)\right\rceil+\binom{r}{2}-1\mbox{.}$$
Let $A=\left(\binom{r}{2}-1\right)\left(\frac{\binom{n}{2}}{g_r(m,n)}-1\right)$.

(i). It is sufficient to show $A>1$ when $n\ge r+1$.
We first prove that $A>1$ when $n\ge 2r-2$. Let $\alpha_r(n)= \frac{g_{r}(0,n)}{\binom{n}{2}}=1-\frac{n-(r-1)}{(r-1)(n-1)}$. Obviously, $\alpha_r(n)$ is decreasing on $n\ge 2r-2$. Hence, $\alpha_r(n)\le\alpha_r(2r-2)=1-\frac{1}{2r-3}$. Therefore,
\begin{eqnarray*}
A &=&\left(\binom{r}{2}-1\right)\left(\frac{\binom{n}{2}}{g_r(m,n)}-1\right)\\
 &\ge&\left(\binom{r}{2}-1\right)\left(\frac{1}{\alpha_r(n)}-1\right)\\
&\ge& \left(\binom{r}{2}-1\right)\left(\frac{1}{1-\frac{1}{2r-3}}-1\right)\\
&=& \frac{r+1}{4}>1\mbox{,}
\end{eqnarray*}
the last inequality holds as $r\ge 4$. Now suppose $r+1\le n\le 2r-3$. Then $k_0=1$, $m\ge 2$, $n=r-1+m$, and $g_r(m,n)=\binom{n}{2}-m=\binom{r-1+m}{2}-m$. Let $h_r(m)=\frac{2g_r(m,r-1+m)}{m}=m+\frac{(r-1)(r-2)}{m}+2r-5$. Then $h_r(m)$ is  a decreasing function of $m$  on $[2,r-2]$. Thus $h_{r}(m)\le h_{r}(2)=2+\frac{(r-1)(r-2)}{2}+2r-5=\frac{(r^2+r-4)}{2}$. Therefore,
\begin{eqnarray*}
A &=& \left(\binom{r}{2}-1\right)\left(\frac{\binom{n}{2}-g_r(m,n)}{g_r(m,n)}\right)\\
&=&\left(\binom{r}{2}-1\right)\cdot\frac{m}{g_r(m,n)}\\
&=&\frac{r^2-r-2}{h_r(m)}\\
&\ge &2\cdot\frac{r^2-r-2}{r^2+r-4}>1\mbox{.}
\end{eqnarray*}

(ii). Since $s\le\max\{r-1,n-1\}$, we can obtain a $T_{r-1}(n-s)$ by deleting at most one vertex in each partition set of $T_{r-1}(n)$. Now let $S$ be the set of the $s$ deleted vertices. It is easy to see that $t_{r-1}(n)-t_{r-1}(n-s)=|G[S]|+|G[S,S^c]|$, where $G=T_{r-1}(n)$. Note that $G[S]=\binom{s}{2}$ since each pair in $S$ comes from different partition sets. To prove (ii), we only need to prove $|G[S,S^c]|\ge\frac{r-2}{r-1}s(n-s)$.

If $s\le m$, then every deleted vertex $v\in S$ comes from a partition set with $k_0+1$ vertices. Hence $|G[\{v\},S^c]|=(n-s)-(k_0+1-1)=n-s-k_0$.  Therefore,
$$|G[S,S^c]|=|S|(n-s-k_0)=s(n-s-\frac{n-m}{r-1})\ge s(n-s-\frac{n-s}{r-1})\ge\frac{r-2}{r-1}s(n-s)\mbox{.}$$
We are done.

If $s>m$, then there are $m$ vertices coming from partition sets with size $k_0+1$ and $s-m$ vertices coming from some partition sets with size $k_0$. Similarly, we have,
\begin{eqnarray*}
|G[S,S^c]| &=& (n-s-k_0)m+(n-s-k_0+1)(s-m)\\
&=&(n-s-\frac{n-m}{r-1})s+s-m\\
&=&\left(n-s-\frac{(n-s)+(s-m)}{r-1}\right)s+s-m\\
&=&\frac{r-2}{r-1}s(n-s)+\left(1-\frac{s}{r-1}\right)(s-m)\\
&\ge&\frac{r-2}{r-1}s(n-s)\mbox{.}
\end{eqnarray*} 
We are done, too.

(iii). We only need to prove the case when $s=1$ since we can write $\frac{t_{r-1}(n-s)}{t_{r-1}(n)}=\prod\limits_{i=n-s+1}^{n}\frac{t_{r-1}(i-1)}{t_{r-1}(i)}$ and $\binom{n-s}{2}/\binom{n}{2}=\prod\limits_{i=n-s+1}^{n}\binom{i-1}{2}/\binom{i}{2}$. Note that when $s=1$, $\binom{n-1}{2}/\binom{n}{2}=1-\frac{2}{n}$, while
$$\frac{t_{r-1}(n-1)}{t_{r-1}(n)}=1-\frac{|T_{r-1}(n)|-|T_{r-1}(n-1)|}{t_{r-1}(n)}=1-\frac{\delta(T_{r-1}(n))}{\frac{1}{2}\sum_{v\in V(T_{r-1}(n))}d(v)}\ge1-\frac{2}{n}\mbox{.}$$
Therefore, (iii) holds.

(iv). When $s=1$, we need to show 
$$\binom{r}{2}-1\ge\frac{r-2}{r-1}\binom{r}{2}+1\mbox{,}$$
this is obviously true since $r\ge 4$.

Now suppose $s\ge 2$. Since $r\ge 4$,  $s\le r-1$ and  $n-s\ge n-r+1\ge 2$, we have $(r-3)(n-s)\ge r-2\ge s-1$. Therefore,
\begin{eqnarray*}
&&(r-2)(n-s)\ge s-1+\frac{2}{s}(n-s)\\
&\Leftrightarrow &(r-2)s(n-s) \ge s(s-1)+2(n-s)\\
&\Leftrightarrow &\left(\frac{1}{2}r-1\right)s(n-s) \ge \binom{s}{2}+(n-s)\\
&\Leftrightarrow &\left(\left(\binom{r}{2}-1\right)-\frac{r-2}{r-1}\binom{r}{2}\right)s(n-s) \ge \left(\binom{r}{2}-\left(\binom{r}{2}-1\right)\right)\binom{s}{2}+n-s\\
&\Leftrightarrow&\left(\binom{r}{2}-1\right)\left(\binom{s}{2}+s(n-s)\right)\ge\binom{r}{2}\left(\binom{s}{2}+\frac{r-2}{r-1}s(n-s)\right)+n-s\\
&\Leftrightarrow&\left(\binom{r}{2}-1\right)\left(\binom{n}{2}-\binom{n-s}{2}\right)=\binom{r}{2}\left(\binom{s}{2}+\frac{r-2}{r-1}s(n-s)\right)+n-s
\mbox{.}
\end{eqnarray*}
The proof is completed.
\end{proof}

\begin{proof}[Proof of Conjecture~\ref{wightmain} for $r\in\{4,5\}$:] Let $\mathcal{G}=(G_0,f)$ be a $k$-weighted graph satisfying the conditions in Conjecture~\ref{wightmain}. 
If $f(e)\le\binom{r}{2}-1$ for any $e\in G_0$, then $f(G_0)\le(\binom{r}{2}-1)|G_0|=(\binom{r}{2}-1)\binom{n}{2}$ and we are done.

Now suppose there exists some edge $e\in G_0$ with $f(e)\ge\binom{r}{2}$. We may suppose that $\mathcal{G}=(G_0,f)$ is a minimum counterexample of Conjecture~\ref{wightmain} with respect to $n$. 
Note that $k\in\{k_1, k_2\}$, where $k_1=\left\lceil\frac{\left(\binom{r}{2}-1\right)\binom{n}{2}}{t_{r-1}(n)}\right\rceil-1$ and $k_2=\left\lceil\frac{\left(\binom{r}{2}-1\right)\binom{n}{2}}{t_{r-1}(n)}\right\rceil=k_1+1$. 

We give more definitions and notation used in the proof. For graphs $G$ and $F$, an {\it $F$-packing} of $G$ is a union of vertex-disjoint copies of $F$ in $G$. The {\em size} of an $F$-packing means the number of disjoint copies of $F$ in this packing. 
For $r\ge 3$ and $2\le s\le r-1$, let $a_{r,s}=\binom{r}{2}-\binom{s+1}{2}+2$. Set $a_{r,r}=1$. Note that $a_{r,r-1}=2$, $a_{r,s-1}-a_{r,s}=s$ and $\binom{r}{2}-a_{r,s-1}+2=\binom{s}{2}$.

First, we pick a maximal $K_{r-1}$-packing $M_{r-1}$ in $G_0$ such that every member $K_{r-1}$ in $M_{r-1}$ has weight sequence bound {${\bf b_{r,r-1}}=(a_{r,r-1}, a_{r,r-2}, a_{r,r-2}+1, \ldots,\binom{r}{2})$.} 
Set $G_{0}^{(r-1)}=G_0-V(M_{r-1})$ and choose a maximal $K_{r-2}$-packing $M_{r-2}$ in $G_0^{(r-1)}$ such that every member $K_{r-2}$ in $M_{r-2}$ has weight sequence bound {${\bf b_{r,r-2}}=(a_{r,r-2}, a_{r,r-3}, a_{r,r-3}+1, \ldots,\binom{r}{2})$.} 
Generally, at the $s$-th step for $2\le s\le r-3$, we pick a maximal $K_{r-s}$-packing $M_{r-s}$ in $G_{0}^{(r-s+1)}$ such that every member $K_{r-s}$ in $M_{r-s}$ has weight sequence bound ${\bf b_{r,r-s}}=(a_{r,r-s}, a_{r,r-(s+1)}, a_{r,r-(s+1)}+1,\ldots, \binom{r}{2})$. Finally, at the $(r-2)$-th step,  we take a maximal  $K_{2}$-packing $M_2$ in $G_{0}^{(3)}$ such that $K_{2}$ has weight sequence bound {${\bf b_{r,2}}=\left(\binom{r}{2}\right)$.}
Let $M_1=G_0-\cup_{i=2}^{r-1} V(M_i)$. Set $m_s=|V(M_s)|$ for $1\le s\le r-1$.
Note that $m_2+m_3+...+m_{r-1}\neq 0$ since there exist edges $e\in G_0$ with $f(e)\ge\binom{r}{2}$.
 Let ${\bf b_{r,r}}=(1, 2,\ldots, \binom{r}{2})$. Then there is no $K_r$ in $G_0$ with weight sequence bound ${\bf b_{r,r}}$ by the assumption of Conjecture~\ref{wightmain}.


Note that our proof is based on induction. Therefore, we need to compute the weights lost when we delete a copy of $K_s$ in $M_s$ from $G_0$. We will address this in the following claims.  

\begin{claim}\label{vertex}
 Let $2\le s\le r-1$. For any vertex $v\in\cup_{i=1}^{s}V(M_i)$ and an $F\cong K_s$ in $M_s$ with $v\notin V(F)$, we have $f(K[\{v\},V(F)])\le (s-1)k+a_{r,s+1}-\delta$ for $k\ge\binom{r}{2}+\delta$, where $\delta\in\{0,1\}$.
\end{claim}
\begin{proof}
We first claim that $K[\{v\}\cup V(F)]$ is not isomorphic to $K_{s+1}$ with weight sequence bound {${\bf b_{r,s+1}}$.} Otherwise, when $s=r-1$, we have a $K_{r}$ in $G_0$ with weight sequence bound ${\bf b_{r,r}}$, a contradiction to our assumption. If $s\le r-2$, then we have another $K_{s+1}$ with weight sequence bound ${\bf b_{r,s+1}}$ disjoint with $\bigcup_{i=s+1}^{r-1}V(M_i)$. Thus we may select a larger $M_{s+1}$. This is a contradiction to the maximality of $M_{s+1}$.
Therefore, $K[\{v\},V(F)]$ should not have weight sequence bound ${\bf c_{r,s}}$, {where ${\bf c_{r,s}}=(a_{r,s+1}, a_{r,s}+1, a_{r,s}+2, ..., a_{r,s-1}-1)$ for $3\le s \le r-1$ and ${\bf c_{r,2}}=(a_{r,3}, a_{r,2})=(\binom{r}{2}-4, \binom{r}{2}-1)$, i.e. ${\bf c_{r,s}}$ is the sequence of integers missed by ${\bf b_{r,s}}$ from ${\bf b_{r,s+1}}$.} Recall that $a_{r,r}=1$.

Suppose the weight sequence of $K[\{v\},V(F)]$ is ${\bf x}=(x_1, x_2, \ldots, x_s)$. Note that $x_i\le k$ for any $1\le i\le s$. 
Denote by ${\bf y}=(y_1, y_2,\ldots, y_s)={\bf c_{r,s}}$. By the above claim, there exists some $j\in [s]$ with $x_j\le y_j-1$. If $j=1$ then $x_1\le y_1-1=a_{r,s+1}-1$. Thus $$f(K[\{v\},V(F)])=\sum_{i=1}^s x_i\le a_{r,s+1}-1+(s-1)k,$$ we are done. Now suppose $2\le j\le s$.
If $s=2$, then  $j=2$ and ${\bf y}={\bf c_{r,2}}=\left(\binom{r}{2}-4, \binom{r}{2}-1\right)$. Thus 
$$f(K[\{v\},V(F)])\le 2(y_2-1)=2\left(\binom{r}{2}-2\right)=\binom{r}{2}+a_{r,3}\le k+a_{r,3}-\delta,$$ the last inequality holds because $\binom{r}{2}\le k-\delta$. Therefore, we are done for $s=2$.
Now suppose $s\ge 3$. Then,  for $1\le i\le j\, (j\ge 2)$, we have $x_i\le x_j\le y_j-1=a_{r,s}+j-2$. Hence
$$f(K[\{v\},V(F)])=\sum_{i=1}^s x_i\le j(a_{r,s}+j-2)+(s-j)k=(s-1)k+j(a_{r,s}+j-2)-k(j-1)\mbox{.}$$
Therefore, it is sufficient to show that $j(a_{r,s}+j-2)-k(j-1)\le a_{r,s+1}-\delta$ in the rest of the proof.

If $s=r-1$, then $a_{r,s+1}=1$ and $a_{r,s}=2$. Since $2\le j\le s=r-1$ and $k\ge\binom{r}{2}\ge r+1\ge j+2$, we have 
$$j(a_{r,s}+j-2)-k(j-1)=j^2-k(j-1)\le j^2-(j+2)(j-1)=2-j\le 0= a_{r,s+1}-1\mbox{.}$$
If $3\le s\le r-2$, then $a_{r,s}=a_{r,s+1}+s+1$. Let $g(x)=(x-1)^2-\left(\binom{s+1}{2}-2+\delta\right)(x-1)$ be a function on $[2,s]$. Since $g$ is convex on $[2,s]$, we have 
$$g(x)\le\max\left\{g(2),g(s)\right\}=\max\left\{3-\delta-\binom{s+1}{2}, 1-\delta-(s-1)\binom{s}{2}\right\}\mbox{.}$$
Since $s\ge 3$, we have $3-\delta-\binom{s+1}{2}\le -(s+\delta)$ and $1-\delta-(s-1)\binom{s}{2}\le -(s+\delta)$. Thus, $g(x)\le -(s+\delta)$. Therefore, 
\begin{eqnarray*}
j(a_{r,s}+j-2)-k(j-1)&=&a_{r,s+1}+s+j-1+(j-1)(a_{r,s}+j-2)-k(j-1)\\
&=&a_{r,s+1}+(j-1)(a_{r,s}+j-k-1)+s\\
&=&a_{r,s+1}-(j-1)\left(k-\binom{r}{2}+\binom{s+1}{2}-1-j\right)+s\\
&\le&a_{r,s+1}-(j-1)\left(\delta+\binom{s+1}{2}-1-j\right)+s\\
&=&a_{r,s+1}+(j-1)^2-\left(\binom{s+1}{2}-2+\delta\right)(j-1)+s\\
&=&a_{r,s+1}+g(j)+s\\
&\le&a_{r,s+1}-(s+\delta)+s\\
&=&a_{r,s+1}-\delta\mbox{.}
\end{eqnarray*}
This completes the proof of Claim~\ref{vertex}.
\end{proof}

Furthermore, we have the following claim when $r=4,5$.
\begin{claim}\label{45}
For $r\in\{4,5\}$, $2\le s\le r-1$, we have $(s-1)k+a_{r,s+1}-1\le\frac{r-2}{r-1}sk$.
\end{claim}
\begin{proof}
It is sufficient to show $a_{r,s+1}-1\le (1-\frac{s}{r-1})k$.
If $s=r-1$, then $a_{r,s+1}=a_{r,r}=1$. Thus $a_{r,s+1}-1=0=(1-\frac{s}{r-1})k$ and we are done.
Else, we have $s+3\ge 5\ge r$. 
Since $k\ge\binom{r}{2}$, we have  $a_{r,s+1}-1=\binom{r}{2}-\binom{s+2}{2}+1=\binom{r}{2}-\frac{s(s+3)}{2}\le\binom{r}{2}-\frac{sr}{2}=(1-\frac{s}{r-1})\binom{r}{2}\le(1-\frac{s}{r-1})k$.
\end{proof}
\begin{claim}\label{n=r}
	$n\ge r+1$.
\end{claim}
\begin{proof}
Note that when $n=r$, $k_2=\binom{r}{2}$, and $k_1=\binom{r}{2}-1<\binom{r}{2}\le k$. Thus, $k$ must be $k_2$ and $\delta=0$ if $n=r$. 
In the following, we show that $f(G_0)\le k\cdot t_{r-1}(n)=k\left({r\choose 2}-1\right)$, that is a contradiction to the fact that $\mathcal{G}=\{G_0, f\}$ is a counterexample. Let ${\bf x}=\left(x_1, x_2, ..., x_{\binom{r}{2}}\right)$ be the weight sequence of $G_0$. Then, there must exist some $j\in[\binom{r}{2}]$ with $x_j<j$. Thus, we get $x_i\le x_j\le j-1$ for $i\le j$ and $x_i\le k_2=\binom{r}{2}$ for $j<i\le\binom{r}{2}$. Therefore,
$$f(G_0)=\sum_{i=1}^{\binom{r}{2}}x_i\le j(j-1)+\left(\binom{r}{2}-j\right)\binom{r}{2}={f_0(j)}\mbox{.}$$
Since $f_0(j)$ is a  convex function on $j\in[\binom{r}{2}]$, we have 
$$f(G_0)\le f_0(j)\le\max\left\{f_0(1), f_0\left(\binom{r}{2}\right)\right\}=\left(\binom{r}{2}-1\right)\binom{r}{2}=k\cdot t_r(n).$$ 
\end{proof}
Now suppose $n\ge r+1$.  By Proposition~\ref{bound} (i), $k_2\ge\binom{r}{2}+1$. We need more computation.
\begin{claim}\label{induction}
For $r\in\{4,5\}$, $\delta\in\{0,1\}$, let $1\le s\le r-1$ be the minimal index with $m_s>0$. Pick an $F\cong K_s$ in $M_s$. If $k\ge\binom{r}{2}+\delta$, then
$$f(F)+f(K[V(F),V(F)^c])\le k\left(\binom{s}{2}+\frac{r-2}{r-1}s(n-s)\right)+(1-\delta)(n-s)\mbox{,}$$
where $V(F)^c=V(G_0)\setminus V(F)$. 
\end{claim}
\begin{proof}
Since $F\cong K_s$, we have $f(F)\le k\binom{s}{2}$.
For any $H\cong K_j$ in $M_j$ with $V(H)\in V(F)^c$, by the minimality of $s$, $s\le j\le r-1$. Let $m_s'=m_s-1\ge 0$ and $m_j'=m_j$ for $j>s$. Then $\sum_{j=s}^{r-1}jm_{j}'=n-s$. By Claims~\ref{vertex} and~\ref{45}, for any $v\in V(F)$, 
$$f(K[\{v\},V(H)])\le (j-1)k+a_{r,j+1}-\delta\le \frac{r-2}{r-1}jk+(1-\delta).$$ 
Thus 
\begin{eqnarray*}
f(K[\{v\},V(F)^c])&=&\sum_{j=s}^{r-1}\sum_{\substack{H\in M_j\\ H\neq F}}f(K[\{v\},V(H)])\\
&\le&\sum_{j=s}^{r-1}m_j'\left(\frac{r-2}{r-1}jk+(1-\delta)\right)\\
&=&\frac{r-2}{r-1}k\sum_{j=s}^{r-1}jm_j'+(1-\delta)\sum_{j=s}^{r-1}m_{j}'\\
&\le&\frac{r-2}{r-1}k(n-s)+\frac{(1-\delta)(n-s)}s\mbox{,}
\end{eqnarray*}
the last inequality holds since $\sum_{j=s}^{r-1}m_{j}'\le\frac{1}{s}\sum_{j=s}^{r-1}jm_j'=\frac{n-s}{s}$.
Therefore, 
\begin{eqnarray*}
f(F)+f(K[V(F),V(F)^c])&\le&k\binom{s}{2}+\sum_{v\in V(F)}f(K[\{v\},V(F)^c])\\
&\le& k\binom{s}{2}+\frac{r-2}{r-1}ks(n-s)+(1-\delta)(n-s)\mbox{.}
\end{eqnarray*}. 
\end{proof}
Now we are ready to prove Conjecture~\ref{wightmain} for $r=4,5$.
\begin{claim}
	Conjecture~\ref{wightmain} holds for $r=4,5$.
\end{claim}	
\begin{proof}
	By Claim~\ref{n=r}, we have $n\ge r+1$. By Proposition~\ref{bound} (i),  $k_2\ge\binom{r}{2}+1$. Recall that the $k$-weighted graph $\mathcal{G}=(G_0,f)$ is a minimum counterexample of Conjecture~\ref{wightmain} and $k\in\{k_1,k_2\}$.  Let $1\le s\le r-1$ be the minimum index with $m_s>0$ and $F\in M_s$ be a copy of $K_s$ in $G_0$  with weight sequence bound ${\bf b_{r,s}}$.

If $k=k_2\ge\binom{r}{2}+1$, then $f(G_0)>k_2\cdot t_{r-1}(n)$ by our assumption that $\mathcal{G}=(G_0,f)$ is a counterexample. Therefore, by Claim~\ref{induction} with $\delta=1$ and Proposition~\ref{bound} (ii),
\begin{eqnarray*}
f(G_0[V(F)^c])&=&f(G_0)-\left(f(F)+f\left(K[V(F),V(F)^c]\right)\right)\\
&>&k_2\cdot\left(t_{r-1}(n)-\left(\binom{s}{2}+\frac{r-2}{r-1}s(n-s)\right)\right)\\
&\ge& k_2\cdot t_{r-1}(n-s)\mbox{.}
\end{eqnarray*} 
Apparently, the restriction of $\mathcal{G}$ on $G_0[V(F)^c]$ is a $k$-weighted graph without any $K_r$ of weight sequence bound $\left(1,2,...,\binom{r}{2}\right)$. However, $G_0[V(F)^c]$ has a smaller number of vertices than $G_0$  while $f(G_0[V(F)^c])>k_2\cdot t_{r-1}(n-s)$, which contradicts to the minimality of $G_0$. 

If $k=k_1$, then we have $f(G_0)>(\binom{r}{2}-1)\binom{n}{2}$ because  $\mathcal{G}=(G_0,f)$ is a counterexample. 
If $k_1=\binom{r}{2}$, applying Claim~\ref{induction} with $\delta=0$ and by Proposition~\ref{bound} (iv), we have
\begin{eqnarray*}
f(G_0[V(F)^c])&=&f(G_0)-(f(F)+f(K[V(F),V(F)^c]))\\
&>&\left(\binom{r}{2}-1\right)\binom{n}{2}-\left(\binom{r}{2}\left(\binom{s}{2}+\frac{r-2}{r-1}s(n-s)\right)+n-s\right)\\
&\ge&\left(\binom{r}{2}-1\right)\binom{n-s}{2}\mbox{.}
\end{eqnarray*}
Clearly,  $G_0[V(F)^c]$ is a smaller counterexample for the same reason as the above case, a contradiction again.
Now suppose $k_1\ge\binom{r}{2}+1$. By Proposition~\ref{bound} (iii), we have $\frac{t_{r-1}(n-s)}{t_{r-1}(n)}\ge\frac{\binom{n-s}{2}}{\binom{n}{2}}$. Hence,
$$\frac{\binom{n}{2}-\binom{n-s}{2}}{t_{r-1}(n)-t_{r-1}(n-s)}\ge\frac{\binom{n}{2}}{t_{r-1}(n)}\mbox{.}$$
By Proposition~\ref{bound} (ii), we have $t_{r-1}(n)-t_{r-1}(n-s)\ge\binom{s}{2}+\frac{r-2}{r-1}s(n-s)$. Then
$$\binom{n}{2}-\binom{n-s}{2}\ge\frac{\binom{n}{2}}{t_{r-1}(n)}\left(\binom{s}{2}+\frac{r-2}{r-1}s(n-s)\right)\mbox{.}$$
Since $k_1=\left\lceil\frac{(\binom{r}{2}-1)\binom{n}{2}}{t_{r-1}(n)}\right\rceil-1<(\binom{r}{2}-1)\frac{\binom{n}{2}}{t_{r-1}(n)}$, we have
$$\left(\binom{r}{2}-1\right)\left(\binom{n}{2}-\binom{n-s}{2}\right)>k_1\left(\binom{s}{2}+\frac{r-2}{r-1}s(n-s)\right)\mbox{.}$$
By Claim~\ref{induction},
\begin{eqnarray*}
f(G_0[V(F)^c])&=&f(G_0)-(f(F)+f(K[V(F),V(F)^c]))\\
&>&\left(\binom{r}{2}-1\right)\binom{n}{2}-k_1\left(\binom{s}{2}+\frac{r-2}{r-1}s(n-s)\right)\\
&\ge&\left(\binom{r}{2}-1\right)\binom{n-s}{2}\mbox{,}
\end{eqnarray*}
a contradiction, too.
\end{proof}
\end{proof}
\section{Discussion and Remarks}
We restate Conjectures~\ref{CONJ: C2} and~\ref{wightmain} in the following. 
 \begin{conj}[Conjecture~1.3]
	Let $n\ge r-1\ge 2$, $k\ge\binom{r}{2}$. 
	Then 
	$${\ex_k(n, K_r)}\le\max\left\{\left(\binom{r}{2}-1\right)\binom{n}{2}, k\cdot t_{r-1}(n)\right\}\mbox{.}$$
\end{conj}

\begin{conj}[Conjecture 2.2]
	Let $n\ge r\ge 3$ and $k\ge\binom{r}{2}$. Let $\mathcal{G}=(G_0,f)$ be a $k$-weighted graph, where $G_0\cong K_n$ on vertex set $[n]$. 
	Suppose that $G_0$ contains no complete subgraph $Kr$ with weight sequence bound $\left(1,2,\ldots,\binom{r}{2}\right)$. Let $k_2=\left\lceil\frac{\left(\binom{r}{2}-1\right)\binom{n}{2}}{t_{r-1}(n)}\right\rceil$ and $k_1=k_2-1$. Then the following hold.
	
	(i) If $k=k_1$, then $f(G_0)\le(\binom{r}{2}-1)\binom{n}{2}$.
	
	(ii) If $k=k_2$, then $f(G_0)\le k\cdot t_{r-1}(n)$.
\end{conj}
In this article, we confirm Conjecture~\ref{CONJ: C2} is true for $r=4,5$ by showing Conjecture~\ref{wightmain} holds for $r\in\{4,5\}$. Our computation will stuck when $r\ge 6$. It will be very interesting to verify Conjecture~\ref{wightmain} for all $r\ge 3$.


\begin{thebibliography}{99}
\bibitem{Caen94}D. de Caen, The current status of Tur\'an's problem on hypergraphs, Extremal
Problems for Finite Sets, Visegr\'ad, 1991, Bolyai Soc. Math. Stud., Vol. 3, pp.
187-197, J\'anos Bolyai Math. Soc., Budapest, 1994.	

\bibitem{K3}
P. Frankl, Graphs without rainbow triangles, arXiv: 2203.07768, 2022.

	

\bibitem{Furedi91} Z. F\"uredi, Tur\'an type problems, Surveys in combinatorics, London Math. Soc.
Lecture Note Ser. 166, Cambridge Univ. Press, Cambridge, 1991, 253-300.
	
\bibitem{GHLSTVZh}E. Gy\"ori, Z. He, Z. Lv, N. Salia, C. Tompkins, K. Varga, and X. Zhu,	Some remarks on graphs without rainbow
triangles, 	arXiv:2204.07567, 2022. 

	


\bibitem{Kee2011}
P. Keevash, Hypergraph Tur\'an problems, Surveys in Combinatorics, Cambridge University Press, 2011 , pp. 83-140.


\bibitem{KSSV04}
P. Keevash, M. Saks, B. Sudakov, J. Verstra\"ete, Multicoloured Tur\'an problems, Advances in Applied Mathematics, 33(2004), 238-262.


\bibitem{sido95}
A. Sidorenko, What we know and what we do not know about Tur\'an numbers,
Graphs and Combinatorics 11 (1995), 179-199.

\bibitem{S68} M. Simonovits, A method for solving extremal problems in graph theory, stability problems, in: Theory of
Graphs, Proc. Colloq., Tihany, 1966, Academic Press, New York, Akad. Kiad\'o, Budapest, 1968, pp. 279-319.

\bibitem{Tr}
P. Tur\'an, On an extremal problem in graph theory, Matematikai \'es Fizikai Lapok (in Hungarian), 48(1941), 436-452.



\end{thebibliography}
\end{document}